# Quadrature rules for $C^0$ and $C^1$ splines, a recipe

Helmut Ruhland[*]

ABSTRACT. Closed formulae for Gaussian or optimal, 1-parameter quadrature rules in a compact interval [a, b] with non uniform, asymmetric subintervals, arbitrary number of nodes per subinterval for the spline classes $S_{2N, 0}$ and $S_{2N+1, 1}$, i.e. even and odd degree are presented.
Also the rules for the 2 missing spline classes $S_{2N-1, 0}$ and $S_{2N, 1}$ (the so called "1/2-rules"), i.e. odd and even degree are presented.
These quadrature rules are explicit in the sense, that they compute the nodes and their weights in the first/last boundary subinterval and, via a recursion the other nodes/weights, parsing from the first/last subinterval to the "middle" of the interval. These closed formulae are based on the semi-classical Jacobi type orthogonal polynomials.

## 1. Introduction

The aim of this paper is to present as a recipe the formulae and the algorithm to create quadrature rules for certain spline classes. Because this paper should serve just as recipe, I tried to make it as short and self contained as possible. For a detailed introduction about this topic see e.g. [5]. In order to the shortness of the paper I did not add any proofs.

In the sequel, we will denote by S be the number of subintervals in [a, b] and by D the degree of the quadrature.

A short description of the problem:

Determine N nodes $x_i$ in the interval [a, b], every node $x_i$ has an assigned weight $w_i$, now the integration of a function f (x) can be approximated by a weighted summation of f (x), evaluated at the nodes $x_i$.

$$\int_a^b f(x)\,dx = \left(\sum_{i=1}^N f(x_i)\,w_i\right) + R(f)$$

It is required, that for f (x) in a certain linear function space, the summation is exact, i.e. R (f) = 0. The function space is here the function space $S_{D, c}$ of splines with degree D and continuity class c. This means in the each of the S subintervals in [a, b] f (x) is a polynomial of maximal degree D, at the inner boundaries of the subintervals the derivations up to the $c^{th}$ derivation are continuous.

The N pairs $(x_i, w_i)$ we call quadrature rules and we are looking for rules with the minimal number of nodes N.

To the best of my knowledge closed formulae are known only for some special cases of the problem. See [5] for the case $C^1$, D = 5, and uniform subintervals, [4] for the case $C^1$, D = 3, and non uniform subintervals, and [3] for the case $C^0$, D = 4, and uniform subintervals.

---

[*] Helmut.Ruhland50@web.de

In the case of the "1/2-rules" see [11] for the case $C^1$, D = 4, and 2 uniform subintervals.

Here we present a set of closed formulae for the case $C^0$ and even D (that includes the case treated by Nikolov in [3]) and a set of formulae for the case $C^1$ and odd D (that includes the cases treated by Bartoň et al in [4] and [5]) for arbitrary even/odd D. In the case $C^0$ the quadrature rules are not Gaussian, but optimal in the sense of the minimal number of nodes and, therefore, we build a 1-parameter family. These formulae are based on semi-classical Jacobi type orthogonal polynomials and it is easy to see that for the "limit" S → 1, i.e. just 1 subinterval, the formulae are those corresponding to the classical Gauss-Legendre quadrature rules.

The 2 cases $C^0$ and odd D, i.e. $S_{2N-1, 0}$, and $C^1$ and even D, i.e. $S_{2N, 1}$, are a little bit more involved. They are very interesting for the case $C^1$, because the Galois group is bigger than the default one (being the direct product of symmetric groups). In such a case, the Galois group is the product of the default group with an additional Abelian 2-group of order $2^{\text{floor }(S/2)}$ as factor. The automorphisms of this group "mix" the zeros of different subintervals.

So the problem of closed formulae for the continuity classes 0 and 1 is solved in a satisfactory way? The case $C^2$, even degree D, Gaussian (but not 1-parameter optimal) non uniform but **symmetric** subintervals seems to be attackable (maybe $C^3$ and odd degree D, too).

The rest of this paper is organized as follows. In section 2 (c = 0) and 3 (c = 1) the polynomials $Q_n$ (…, x) and $M_n$ (…, x) are presented (… as arguments here means one or two arguments, depending on c). These polynomials play a similar role as the Legendre polynomials $P_n$ (x) in the classical Gauss-Legendre quadrature, i.e. the nodes are the roots of such a polynomial and the weights are given in the same style as for Gauss-Legendre (the product of a derivation of $Q_n$ (…, x) and $Q_{n-1}$ (…, x) in the denominator of the formula for the weights). In section 4 follows the algorithm, which is essentially the same for both c. To implement the algorithm section 2 or 3 depending on the case c and section 4 should do.
Section 5 contains an example for a non-uniform, asymmetric case to visualize the algorithm and which could also be used to check an implementation. 2 further examples with quadrature rules for the real line are added.

For the "1/2-rules" in section 6 (c = 0) and 7 (c = 1) the additional quantities apart from the quantities in section 2 (c = 0) and 3 (c = 1) are given. In section 8 follows the algorithm, which is essentially the same for both c.

Section 9 contains an example for a non-uniform, asymmetric "1/2-rule" case to visualize the algorithm and which could also be used to check an implementation.

Section 10. shows how the nodes can be distributed in the subintervals. It visualizes the different types of optimal and suboptimal quadrature rules covered by this recipe.

Section 11 has some comments on a proof and the numerical verifications I have done. This section also shows, how $Q_n$ (…, x) and $M_n$ (…, x) are related to the semi-classical Jacobi type polynomials.

## 2. Continuity class c = 0, formulae for $S_{2N, 0}$

### 2.1. Polynomials and weights for the subintervals left/right from the interval $S_{M\ =\ middle}$:

The polynomials $Q_n(\alpha, x)$ w i t h o u t even/odd symmetry
Let $P_n(x)$ be the Jacobi polynomials $P_n^{(1, 0)}(x)$

(1) $\quad \begin{aligned} F(n) &= 1 + \alpha\, n\,(n+1) \\ Q_n &= (F(n) + \alpha\, n)\, P_n + \alpha\,(1 - x)\, D(P_n) \end{aligned}$

For the 1st subinterval i.e. $\alpha = 0$: $F(n) = 1$, then $Q_n(0, x) = P_n(x)$

(2) $\quad w_i = \dfrac{2\,(2n+1)\,F(n)^2}{n\,(n+1)\,(Q_n)'(x_i)\,Q_{n-1}(x_i)\,(1 - x_i)}$

### 2.2. The recursion map $\alpha_{i-1} \to \alpha_i$ to step from a subinterval to the next:

(3) $\quad \begin{aligned} \Gamma &= (n+1)\,(1 + n\,(n+2)\,\alpha) \\ \alpha_{new} &= \dfrac{1 + (n+1)^2\,\alpha}{(n+1)\,\Gamma} \end{aligned}$ a fractional linear map

An additional monomial transformation if a stretching factor $\lambda \neq 1$ is introduced. Indeed, let $\lambda_i = L_i / L_{i-1}$ ($L_i$ denotes the length of the i-th subinterval) be the stretching factor (in the uniform case, $\lambda_i = 1$)

(3.1) replace the $\alpha_{new}$ in (3) above by $\alpha_{new} / \lambda_i$

### 2.3. Polynomial and weights for the subinterval in the "middle" $S_M$ with $N + 1$ nodes:

The polynomials $M_n(\alpha_L, \alpha_R, x)$ w i t h even/odd symmetry for $\alpha_L = \alpha_R$
Let now $P_n(x)$ be the Legendre polynomials i.e. $P_n^{(0, 0)}(x)$

(5) $\quad \begin{aligned} H(n) &= 1 + n^2\,(\alpha_L + \alpha_R + (n-1)(n+1)\,\alpha_L\,\alpha_R) \\ H_1(n) &= \alpha_L + \alpha_R + 2n(n+1)\,\alpha_L\,\alpha_R \\ H_2(n, \alpha) &= 1 + \alpha\, n\,(n+1) \\ M_n &= (H(n) + n\,H_1(n))\,P_n + (\alpha_L\,H_2(n, \alpha_R)(1 - x) - \alpha_R\,H_2(n, \alpha_L)(1 + x))\,D(P_n) \\ M_{\omega, n} &= M_n + \omega\,M_{n-1} \end{aligned}$

For the 1st subinterval (and the only one) i.e. $\alpha_L, \alpha_R = 0$: $H(n) = 1$ and $M_n(0, 0, x) = P_n(x)$ the Legendre polynomials

(6) $\quad w_i = \dfrac{2\,H(n)^2}{n\,(M_{\omega, n})'(x_i)\,M_{n-1}(x_i)}$

## 3. Continuity class c = 1, formulae for $S_{2N+1,\,1}$

### 3.1. Polynomials and weights for the subintervals left/right from the interval $S_{M\,=\,middle}$:

The polynomials $Q_n(\alpha, \beta, x)$ w i t h o u t  even/odd symmetry
Let $P_n(x)$ be the Jacobi polynomials $P_n^{(2,\,0)}(x)$. In this situation

(1)
$$F(n) = 1 + n(n+2)(\alpha + 6(n^2 + 2n - 1)\beta - 3(n-1)n(n+1)^2(n+2)(n+3)\beta^2)$$
$$F_1(n) = \alpha + 12\beta(n^2 + 3n + 1 - n(n+1)^2(n+2)^2(n+3)\beta)$$
$$F_2(n) = \beta(1 - 3n(n+1)(n+2)(n+3)\beta)$$
$$Q_n = (F(n) + n F_1(n))P_n + F_1(n)(1-x)D(P_n) - 36 F_2(n)D(P_n)$$
$$+ 12 F_2(n)(1-x)D(P_n)^2$$

For the 1st subinterval i.e. $\alpha, \beta = 0$: F(n) = 1, $F_1(n)$, $F_2(n) = 0$, then $Q_n(0, 0, x) = P_n(x)$

(2)
$$w_i = \frac{8(n+1)F(n)^2}{n(n+2)(Q_n)'(x_i)Q_{n-1}(x_i)(1-x_i)^2}$$

### 3.2. The recursion map $\alpha_{i-1}, \beta_{i-1} \to \alpha_i, \beta_i$ to step from a subinterval to the next:

(3)
$$\Gamma = (1 + n(n+3)\alpha + 6n(n+3)(n^2 + 3n - 1)\beta$$
$$- 3n^2(n-1)(n+4)(n+2)(1+n)(n+3)^2\beta^2)(1+n)(n+2)/2$$
$$E = 1 + (n+1)(n+2)(\alpha + 3n(n+3)\beta(2 - (n-1)(n+1)(n+2)(n+4)\beta))$$
$$G = 1 - 3n(n+1)(n+2)(n+3)\beta$$
$$\alpha_{new} = -\alpha + E(8n^2 + 24n + 12 + n(n+3)((11n^2 + 33n + 16)\alpha$$
$$+ 12(4n^4 + 24n^3 + 34n^2 - 6n - 8)\beta + 3n(n+1)(n+2)(n+3)($$
$$-4(n+1)(n+2)(2n^2 + 6n - 5)\beta^2$$
$$- 3(n-1)n(n+1)(n+2)(n+3)(n+4)\alpha\beta^2 + 2(3n^2 + 9n - 6)\alpha\beta + \alpha^2)))$$
$$/(12\,\Gamma^2)$$
$$\beta_{new} = \beta + \frac{E\,G}{6(n+1)(n+2)\Gamma}$$

An additional monomial transformation if a stretching factor $\lambda \neq 1$:
Let $\lambda_i = L_i / L_{i-1}$ ($L_i$ the length of the $i^{th}$ subinterval) be the stretching factor

(3.1)      replace the $\alpha_{new}$ in (3) above by $\alpha_{new} / \lambda_i$ and $\beta_{new}$ by $\beta_{new} / \lambda_i^2$

### 3.3. Polynomial and weights for the subinterval $S_M$ in the middle with $N + 1$ nodes:

The polynomials $M_n(\alpha_L, \alpha_R, \beta_L, \beta_R, x)$ with even/odd symmetry for $\alpha_L = \alpha_R$, $\beta_L = \beta_R$
Let now $P_n(x)$ be the Legendre polynomials i.e. $P_n^{(0,0)}(x)$. In this situation

(5)
$$H_0(n, \alpha, \beta) = 1 + n(n-1)(\alpha + (n+1)(n-2)\beta(6 - 3\beta(n+2)n(n-1)(n-3)))$$
$$H1L_0(n) = H_0(n+1, \alpha_L, \beta_L), \quad H1R_0(n) = H_0(n+1, \alpha_R, \beta_R)$$
$$H(n) = \frac{1}{2} H_0(n, \alpha_L, \beta_L) H1R_0(n) + \frac{1}{2} H_0(n, \alpha_R, \beta_R) H1L_0(n)$$
$$\quad - 36(n-1)n^2(n+1)(\beta_L - \beta_R)^2$$
$$H_1(n, \alpha, \beta) = \alpha + 12 n(n+1)\beta(1 - (n-1)(n+2)(n^2 + n + 3)\beta)$$
$$H_2(n, \beta) = \beta(1 - 3(n-1)n(n+1)(n+2)\beta)$$
$$H_3(n, \alpha, \beta) = H_0(n+1, \alpha, \beta) + 24 n(n+1) H_2(n, \beta)$$
$$H_4(n, \alpha, \beta) = 1 + n(n+1)(2\alpha + 3(n-1)n(n+1)(n+2)(13n^2 + 13n - 18)\beta^2)$$

$$\begin{aligned}
M_n &= \frac{1}{2}(H_3(n, \alpha_L, \beta_L) H1R_0(n) + H_3(n, \alpha_R, \beta_R) H1L_0(n)) P_n \\
&+ (H_1(n, \alpha_L, \beta_L) H1R_0(n)(1-x) - H_1(n, \alpha_R, \beta_R) H1L_0(n)(1+x)) D(P_n) \\
&+ 12 (H_2(n, \beta_L) H1R_0(n)(1-x) + H_2(n, \beta_R) H1L_0(n)(1+x)) D(P_n)^2 \\
&- 36 (\beta_L - \beta_R)^2 n(n+1)(n(n+1) P_n - 2x D(P_n) + 2 D(P_n)^2) \\
&+ 12 (H_2(n, \beta_R) H_4(n, \alpha_L, \beta_L) - H_2(n, \beta_L) H_4(n, \alpha_R, \beta_R)) D(P_n) \\
&+ 72 n(n+1)(\beta_L - \beta_R)(\beta_L + \beta_R - 6(n-1)n(n+1)(n+2)\beta_L \beta_R) x D(P_n)^2
\end{aligned}$$

For the 1st subinterval (and the only one) i.e. $\alpha, \beta = 0$: $H(n) = 1$, $H_1(n)$, $H_2(n) = 0$ and, then $M_n(0, 0, 0, 0, x) = P_n(x)$

(6)
$$w_i = \frac{2 H(n)^2}{n (M_n)'(x_i) M_{n-1}(x_i)}$$

# 4. Algorithm to calculate the nodes and weights for $S_{2N,\,0}$, $S_{2N+1,\,1}$

---

| | | |
|---|---|---|
| **Algorithm** | | Gaussian/Optimal Quadrature ([a, b], N, S, $S_M$) spline classes $S_{2N,\,0}$ and $S_{2N+1,\,1}$ |

---

1: **Input:** Compact interval [a, b]
    N = number of nodes per subinterval
    S = number of non uniform subintervals
       a list with the S subinterval lengths
    $S_M$ = the subinterval in the "Middle"
       $1 \leq S_M \leq S$ with N + 1 nodes

---

2:     calculate by (1) the polynomials $Q_{N-1}$, $Q_N$
3:     set the c+1 initial values for the parameters α, … to 0

4:     **for** s from 1 to $S_M$ - 1 **do**
5:         calculate the N roots $x_i$ of $Q_N$, the weights $w_i$ with (2) for each $x_i$
        comment: to scale, the roots $x_i$ of $Q_N$ lie (have to lie) in the interval [-1, +1]
6:         **out:** scale the N 2-tupels (+ $x_i$, $w_i$) to subinterval s
7:         calculate the new c+1 parameters α, … for the next subinterval with the recursion (3)
8:         if stretching, apply the additional monomial transformation (3.1)
9:     **end for**
10:     $α_L$, … = α, …, Index L = left

11:     Repeat steps 2 – 9 for the subintervals right from the $S_M$, i.e. s decreasing
       from S to $S_M$ + 1 now with the reflected 2-tupels (- $x_i$, $w_i$) to subinterval s
12:     $α_R$, … = α, …, Index R = right

14:     calculate by (5) the polynomials $M_N$, $M_{N+1}$ or $M_N$, $M_{ω,N+1}$ (for c = 0, ω a free parameter)
15:     calculate the (N + 1) roots $x_i$ of $M_{N+1}$ or $M_{ω,N+1}$ the weights $w_i$ with (6) for each $x_i$
16:     **out:** scale the (N + 1) 2-tupels ($x_i$, $w_i$) to subinterval [$S_M$]

18:     **Output:** the set { $x_i$, $w_i$ } $_{i\,=\,1\,…\,N\,S\,+\,1}$ of nodes and weights

---

*Remark: For S = number of uniform subintervals = 1, i.e. no inner knots, i.e. no conditions on smoothness, this is of course the classical Gauss-Legendre quadrature, the coefficients H (n) in formulae (6) are 1*

# 5. Some examples

## 5.1. A non uniform, odd degree, optimal, asymmetric $C^1$ quadrature rule

The interval [0, 9] is subdivided in the 6 subintervals [0, 1], [1, 3], [3, 6], [6, 7], [7, 8] and [8, 9] and the degree of the rule D = 3. i.e. the knot vector [ 0,0,0,0, 1,1, 3,3, 6,6, 7,7, 8,8, 9,9,9,9 ]

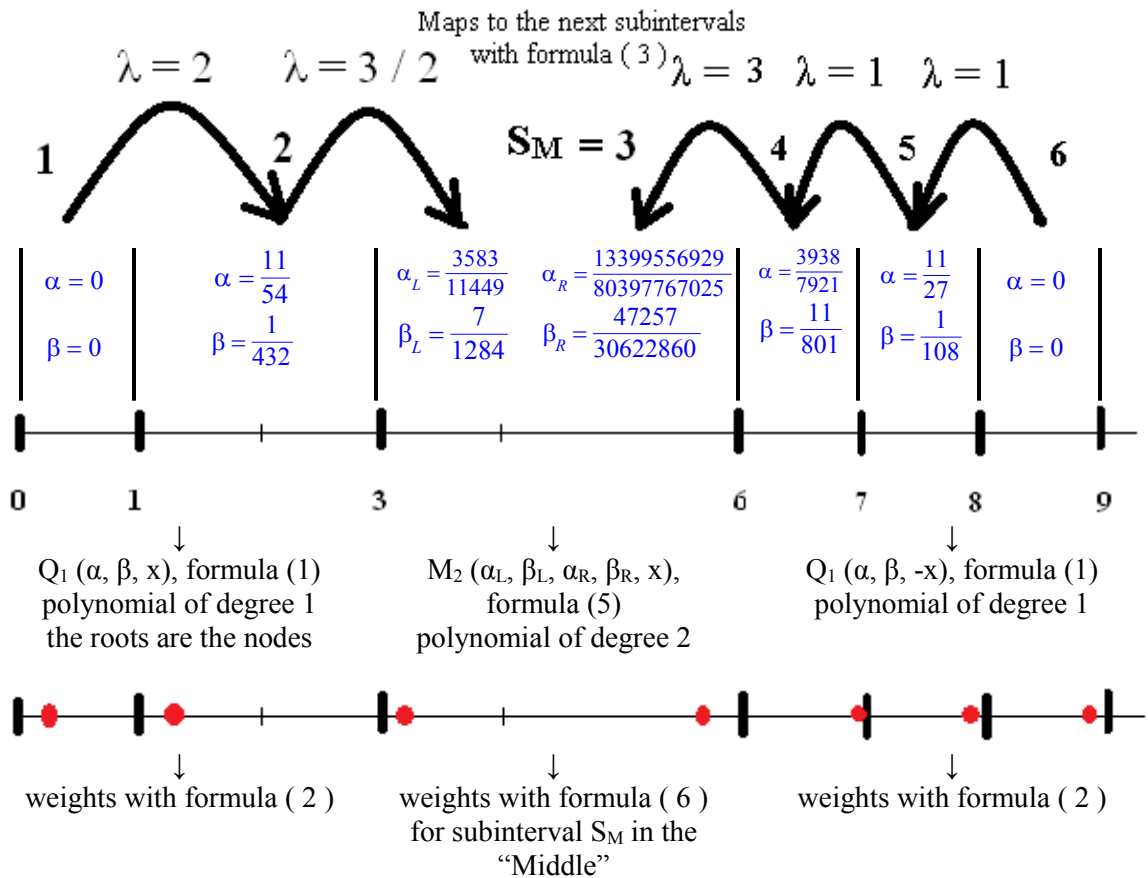

The node/weight list (as rational/algebraic numbers):

(5.1)
```
SQ = sqrt (3556830148073443658426871391555)

# subinterval 1
[ 1 / 4, 16 / 27 ]

# subinterval 2
 [ 76 / 61, 453962 / 309123 ]

# subinterval 3
 [ 922485522061455153 / 210841059447710038 - 135/210841059447710038 * SQ,
 6807081574081001530033959853904 / 9821909391090899005165052208283642169650184625 * SQ
 + 1361950761199921 / 613649356446150]
 [ 922485522061455153 / 210841059447710038 + 135/210841059447710038 * SQ,
 - 6807081574081001530033959853904 / 9821909391090899005165052208283642169650184625 * SQ
 + 1361950761199921 / 613649356446150]

# subinterval 4
 [ 662139 / 94604, 211674482615216 / 212276904201675 ]

# subinterval 5
[ 733 / 92, 194672 / 213867 ]

# subinterval 6
 [ 35 / 4, 16 / 27 ]
```

These results coincide with the result of Johannessen's MATLAB Code, see [9] and [10]

## 5.2. A uniform $C^0$ quadrature rule, even degree, $S \to \infty$ i.e. for the real line

The proof is left for the reader as an exercise.
Hint: Determine the fixed point in the map (3), there is only one attractor with linear convergence, then use formulae (1) and (2). See appendix A.1 with the formula that expresses $Q_n (\alpha, x)$ by Gegenbauer polynomials.

This is the optimal quadrature rule give in [14] in section 2.2 formulae (2.2.1) …

(10) $\quad C_n (x) = C_n^{(3/2)} (x) \quad\quad$ a Gegenbauer polynomial for the weight $(1 - x^2)$

(11) $\quad \delta = \sqrt{\dfrac{n+2}{n}}$

(12) $\quad R_n = C_n(x) + \delta\, C_{n-1}(x)$

(13) $\quad S_n = C_{n-1}(x)\,(2n + 1 + \delta\, x\, n) - C_{n-2}(x)\,(n+1)\,x$

Let the n nodes be:
(14) $\quad x_1 \ldots x_n \quad$ the n roots of the polynomial $R_n (x)$

Let the n weights be:
(15) $\quad w_i = \dfrac{2\,(n+1)\,(2n+1)}{(R_n')(x_i)\, S_n(x_i)} \quad\quad$ for $x_i$ ($i = 1 \ldots n$) a root of $R_n (x)$

In the contrary to the following case c = 1 the nodes are not symmetric (because of the term with $C_{n-1}$ in $R_n$ in (12)). Changing the sign of the square root for δ in (12) is just a reflection at 0 of the roots $x_i$, i.e. $x_i \to -x_i$, because the $C_i (x)$ are even/odd functions

## 5.3. A uniform $C^1$ quadrature rule, odd degree, $S \to \infty$ i.e. for the real line

The proof is left for the reader again as an exercise.
Hint: Determine the rational fixed point with quadratic convergence in the map (3), then use formulae (1) and (2). See appendix A.2 with the formula that expresses $Q_n (\alpha, \beta, x)$ by Gegenbauer polynomials. $Q_n (x)$ has here a factor $(1 + x)$, proof that $A\, Q_n (x) / (1 + x)$ for a constant A to be determined fulfils the 3 term recursion for the Gegenbauer polynomial $C_{n-1} (x)$

There is a second fixed point in the map (3) resulting in further quadrature rule with the same parameters

This is the optimal quadrature rule give in [14] in section 3.1.1. formulae (3.1.1) …

(16) $\quad C_{n-1} (x) = C_{n-1}^{(5/2)} (x) \quad\quad$ a Gegenbauer polynomial for the weight $(1 - x^2)^2$

Let the n nodes be:
(17) $\quad x_1 = -1 \quad\quad x_2 \ldots x_n \;$ the n-1 roots of the Gegenbauer polynomial $C_{n-1} (x)$

Let the n weights be:

(18)
$$w_1 = \frac{16(2n^2 + 6n + 1)}{3n(n+1)(n+2)(n+3)} \qquad \text{for } x_1 = -1$$

$$w_i = \frac{2}{9} \frac{n(n+1)(n+2)}{(C_{n-1}')(x_i) C_{n-2}(x_i)(1-x_i^2)^2} \qquad \text{for } x_i \ (i = 2 \ldots n) \text{ a root of } C_{n-1}(x)$$

These are for low n the well known quadrature rules $S_{2n+1, 1}$ for the real line

- see [5] for the quintic case n = 2 page 15, (39)
- see [6] for the septic case n = 3 page 15, (37), for the nonic case n = 4 page 17, (38)

# 6. Continuity class c = 0, additional formulae for $S_{2N-1, 0}$, the "1/2-rules"

## 6.1. Additional formulae for the subintervals left/right from the $S_M$ ($S_M + 1$):

These formulae are additional to the formulae in section 2. because an iteration step now produces the nodes/weights of 2 subintervals, $\lambda$ is the stretching factor between the 2 subintervals of one iteration step

$$(1.1^*) \quad \omega = -\frac{n(1+(n+1)^2\alpha) + \lambda(n+1)(1+n(n+2)\alpha)}{(n+1)(1+n^2\alpha) + \lambda n(1+(n-1)(n+1)\alpha)}$$

$$(1.2^*) \quad Q_{\omega, n} = Q_n + \omega\, Q_{n-1}$$

## 6.2. The additional extended recursion map $\alpha_{i-1} \to \alpha_i$, depending on $\omega$:

$$\Gamma = (n+1)(1+n(n+2)\alpha) + \omega\, n(1+(n-1)(n+1)\alpha)$$

$$(3.2^*) \quad \alpha_{new} = \frac{n(1+(n+1)^2\alpha) + \omega(n+1)(1+n^2\alpha)}{n(n+1)\,\Gamma}$$

For $\omega = 0$ this is (3) in section 2.

## 6.3. Additional formulae for 2 subintervals $S_M$, $S_M + 1$ in the "middle":

For $S = 0 \bmod 2$ there are 2 subintervals in the "middle", therefore additional $\alpha$ are necessary, $\omega$ is the free parameter because this rule is 1-parameter optimal

$(5.1^*) \quad \alpha_{R, M} = \omega$            for the subinterval $S_M$

$(5.2^*) \quad \alpha_{L, M} = -\dfrac{\alpha_{R, M}}{\lambda}$        for the subinterval $S_M + 1$
                                                 $\lambda$ = stretching factor

# 7. Continuity class c = 1, additional formulae for $S_{2N,\,1}$, the "1/2-rules"

## 7.1. Additional formulae for the subintervals left/right from the $S_M$ ($S_M + 1$):

These formulae are additional to the formulae in section 2. because an iteration step now produces the nodes/weights of 2 subintervals, λ is the stretching factor between the 2 subintervals of one iteration step, here determining ω:

Because the variables A, B, C, depending on α, β, λ and defining a quadratic equation for ω are a little bit long, I do not write it here, see appendix B. The quadratic equation for ω is:

$$(1.1^*) \quad A\omega^2 + B\omega + C = 0$$
$$\omega = \frac{-B + \sqrt{B^2 - 4AC}}{2A}$$

$$(1.2^*) \quad Q_{\omega,\,n} = Q_n + \omega\, Q_{n-1}$$

## 7.2. The additional extended recursion map $\alpha_{i-1}, \beta_{i-1} \to \alpha_i, \beta_i$ depending on ω:

$(3.2^*)$
$$\Gamma = -(n+2)\,(\alpha\,n\,(n+3) - 3\,\beta^2\,n^2\,(n-1)\,(n+4)\,(n+2)\,(n+1)\,(n+3)^2$$
$$+ 6\,\beta\,n\,(n+3)\,(n^2 + 3n - 1) + 1) + \omega\,(-\alpha\,n\,(n+2)\,(n-1) +$$
$$3\,\beta\,n\,(n+2)\,(n-1)\,(\beta\,n\,(n-2)\,(n+3)\,(n+1)\,(n+2)\,(n-1) - 2n^2 - 2n + 6)$$
$$- n)$$

$$A = \omega^2\,(\alpha^2\,n\,(n-1)\,(n+2)\,(n+1)\,(2n^2 + 2n - 3) + \alpha\,($$
$$-3\,\beta^2\,n^2\,(2n^2 + 2n - 9)\,(2n^2 + 2n - 3)\,(n-1)^2\,(n+2)^2\,(n+1)^2$$
$$+ 6\,\beta\,n\,(n-1)\,(n+2)\,(n+1)\,(2n^2 + 2n - 5)\,(2n^2 + 2n - 3)$$
$$+ (2n^2 + 2n - 3)\,(2n^2 + 2n - 1))$$
$$+ 9\,\beta^4\,n^4\,(n-2)\,(n+3)\,(2n^2 + 2n - 9)\,(n+2)^3\,(n-1)^3\,(n+1)^4$$
$$- 36\,\beta^3\,n^2\,(2n^4 + 4n^3 - 9n^2 - 11n + 6)\,(n^2 + n - 3)\,(n-1)^2\,(n+2)^2\,(n+1)^2 +$$
$$6\,\beta^2\,n\,(n-1)\,(n+2)\,(n+1)\,(10n^6 + 30n^5 - 35n^4 - 120n^3 + 67n^2 + 132n - 72)$$
$$+ 12\,\beta\,(n+2)\,(n-1)\,(2n^2 + 2n - 3)\,(n^2 + n - 1) + 2n^2 + 2n - 1) + \omega\,($$
$$2\,\alpha^2\,n\,(n+2)\,(2n^2 + 4n - 3)\,(n+1)^2 + \alpha\,($$
$$-6\,\beta^2\,n^2\,(n+3)\,(n-1)\,(4n^4 + 16n^3 + 20n^2 + 8n - 27)\,(n+2)^2\,(n+1)^2$$
$$+ 12\,\beta\,n\,(n+2)\,(n+1)^2\,(2n^2 + 4n - 3)^2 + 2\,(2n^2 + 4n - 1)\,(2n^2 + 4n + 3))$$
$$+ 18\,\beta^4\,n^3\,(2n^4 + 8n^3 - 5n^2 - 26n + 12)\,(n+3)^2\,(n-1)^2\,(n+2)^3\,(n+1)^4$$
$$- 72\,\beta^3\,n^2\,(n+3)\,(n-1)\,(2n^2 + 4n - 7)\,(n^4 + 4n^3 + 4n^2 - 3)\,(n+2)^2\,(n+1)^2$$
$$+ 12\,\beta^2\,n\,(n-1)\,(n+3)\,(n+2)\,(10n^4 + 40n^3 - n^2 - 82n + 30)\,(n+1)^2$$
$$+ \beta\,(72 + 600\,n^4 - 192\,n + 288\,n^5 + 480\,n^3 + 48\,n^6) + 6 + 8n + 4n^2)$$
$$+ \alpha^2\,n\,(n+3)\,(n+2)\,(n+1)\,(2n^2 + 6n + 1) + \alpha\,($$

$$-3\beta^2 n^2 (2n^2+6n-5)(2n^2+6n+1)(n+3)^2(n+2)^2(n+1)^2$$
$$+6\beta n(n+3)(n+2)(n+1)(2n^2+6n-1)(2n^2+6n+1)$$
$$+(2n^2+6n+1)(2n^2+6n+3))$$
$$+9\beta^4 n^3 (n-1)(n+4)(2n^2+6n-5)(n+3)^3(n+2)^4(n+1)^4$$
$$-36\beta^3 n^2 (2n^4+12n^3+15n^2-9n-8)(n^2+3n-1)(n+3)^2(n+2)^2(n+1)^2$$
$$+$$
$$6\beta^2 n(n+3)(n+2)(n+1)(10n^6+90n^5+265n^4+240n^3-53n^2-24n+12)$$
$$+12\beta n(n+3)(n^2+3n+1)(2n^2+6n+1)+2n^2+6n+3$$
$$B = \omega(-\alpha n(n+2)(n+1)+3\beta^2 n^3(n-1)(n+2)^2(n+1)^3$$
$$-6\beta n(n+2)(n+1)(n^2+n-1)-2-n)-\alpha n(n+2)(n+1)$$
$$+3\beta^2 n^2 (n+3)(n+2)^3(n+1)^3 - 6\beta n(n+2)(n+1)(n^2+3n+1) - n$$

$$\alpha_{new} = \frac{4A}{3(n+1)^2 \Gamma^2}$$

$$\beta_{new} = \frac{B}{3\Gamma(n+1)^2(n+2)n}$$

The $\Gamma$ here is for $\omega = 0$ apart from a factor the $\Gamma$ in (3) in section 3.

### 7.3. Additional formulae for 2 subintervals $S_M$, $S_M + 1$ in the "middle":

For even S there are 2 subintervals in the "middle", therefore additional α, β are necessary. These additional α, β are the $α_{R, M}$, $β_{R, M}$ for the right (the index R) side of the left subinterval and the $α_{L, M}$, $β_{L, M}$ for the left side (the index L) of the right subinterval. These depend on the 6 variables n, λ, $α_L$, $β_L$, $α_R$, $β_R$ including a square root.

Here the defining equations for $α_{R, M}$, $β_{R, M}$

(5.0.1*)
$$A(α, β) = 9 n^2 (n + 2)^2 (n + 1)^4 (-α (n + 3)(n - 1)$$
$$+ 3 β^2 n (n + 4)(n - 2)(n + 2)(n + 3)^2 (n - 1)^2$$
$$- 6 β (n + 3)(n - 1)(n^2 + 2n - 4) - 1)$$

$$B(α, β) = 18 n (n + 2)(n + 1)^2 (α (n^2 + 2n - 1)$$
$$- 3 β^2 n (n - 1)(n + 3)(n + 2)(n^2 + 2n - 4)(n + 1)^2$$
$$+ 6 β (n^2 + 2n - 2)(n^2 + 2n - 1) + 1)$$

$$C(α, β) = 3 α (n + 1)^2 - 9 β^2 n^2 (n + 2)^2 (n + 1)^4 + 18 β n (n + 2)(n + 1)^2 + 3$$

$$DD(α, β) = 3 (n + 1)^2 (α n (n + 2) - 3 β^2 n^2 (n + 3)(n - 1)(n + 2)^2 (n + 1)^2$$
$$+ 6 β n (n + 2)(n^2 + 2n - 1) + 1)$$

$$A_1 = A(α_L, β_L), B_1 = B(α_L, β_L), C_1 = C(α_L, β_L), DD_1 = DD(α_L, β_L)$$

$$A_2 = A(α_R, β_R), B_2 = λ^2 B(α_R, β_R), C_2 = λ^4 C(α_R, β_R), DD_2 = -λ^3 DD(α_R, β_R)$$

(5.0.2*) These $A_i$, $B_i$, $C_i$, $DD_i$ for i = 1, 2 define 2 equations quadratic in $β_{R, M}$ and linear in $α_{R, M}$:

$$A_i β_{R, M}^2 + B_i β_{R, M} + C_i + DD_i α_{R, M} = 0$$

Eliminating $α_{R, M}$ we get the following $A_3$, $B_3$, $C_3$, eliminating $β_{R, M}^2$ we get the following $B_4$, $C_4$, $DD_4$:

(5.0.3*)
$$A_3 = DD_2 A_1 - DD_1 A_2, B_3 = DD_2 B_1 - DD_1 B_2, C_3 = DD_2 C_1 - DD_1 C_2$$
$$B_4 = A_2 B_1 - A_1 B_2, C_4 = A_2 C_1 - A_1 C_2, DD_4 = DD_1 A_2 - DD_2 A_1$$

Replacing now the n, λ, $α_L$, $β_L$, $α_R$, $β_R$ by its rational or numerical values, the equations above get very small and it is no problem to solve the equation system

(5.1*)
$$SQ = B_3^2 - 4 A_3 C_3$$

$$β_{R, M} = \frac{1}{2} \frac{-B_3 - \sqrt{SQ}}{A_3} \quad \text{for the subinterval } S_M$$

$$α_{R, M} = -\frac{C_4 + B_4 β_{R, M}}{DD_4}$$

(5.2*)
$$\alpha_{L,M} = -\frac{\alpha_{R,M}}{\lambda}$$

$$\beta_{L,M} = \frac{\beta_{R,M}}{\lambda^2}$$

for the subinterval $S_M + 1$
$\lambda$ = the stretching factor

# 8. Algorithm to calculate the nodes and weights for $S_{2N-1,\,0}$, $S_{2N,\,1}$ i.e. the "1/2-rules"

| **Algorithm** | Gaussian/Optimal Quadrature ([a, b], N, S, $S_M$) |
|---|---|
| | spline classes $S_{2N-1,\,0}$ and $S_{2N,\,1}$ |
| | **for "1/2-rules"** |

1: **Input:** Compact interval [a, b]
    N = number of nodes per subinterval
    S = number of (non-)uniform subintervals,
        a list with the S subinterval lengths
    $S_M$ = the subinterval in the "middle"
        $1 \leq S_M \leq S$ with N nodes
        $S_M$ has to be odd, if S is even, there
        are 2 subintervals in the "middle"

2: calculate by (1) the polynomials $Q_{N-2}$, $Q_{N-1}$, $Q_N$
3: set the c+1 initial values for the parameters α, … to 0

4: **for** s **from** 1 **to** $S_M - 2$ **by** 2 **do**
    comment: each of these iteration step produces the $x_i$, $w_i$ for the subintervals (s, s+1)

5:     calculate ω with the additional (1.1*), λ = 1 if uniform, else λ = $L_{s+1}$ / $L_s$, L is a length
    (for c = 1 select a sign for the square root in the formula for ω)
6:     calculate the N roots $x_i$ of $Q_{\omega,\,N}$, the weights $w_i$ with (2) for each $x_i$
    comment: to scale, the roots $x_i$ of $Q_{\omega,\,N}$ lie (have to lie) in the interval [-1, +1]
7:     **out:** scale the N 2-tupels (+ $x_i$, $w_i$) to subinterval s
8:     calculate the c+1 parameters α, … for the next subinterval s with (3.2*)
9:     if stretching, apply the additional monomial transformation (3.1)

10:     calculate the N-1 roots $x_i$ of $Q_{N-1}$, the weights $w_i$ with (2) for each $x_i$
    comment: to scale, the roots $x_i$ of $Q_{N-1}$ lie (have to lie) in the interval [-1, +1]
11:     **out:** scale the N-1 2-tupels (+ $x_i$, $w_i$) to subinterval s+1

12:     calculate the new c+1 parameters α, … for the next iteration step with (3)
13:     if stretching, apply the additional monomial transformation (3.1)
14: **end for**
15: for the intervals in the "middle": $α_L$, … = α, …, Index L = left

16: Repeat steps 2 – 14 for the subintervals right from the $S_M$, i.e. s decreasing
    from S to $S_M + 1$ for S odd or $S_M + 2$ for S even now with the reflected 2-tupels (- $x_i$, $w_i$)
17: for the intervals in the "middle": $α_R$, … = α, …, Index R = right

18: **if** S = odd **then**
    comment: just 1 subinterval $S_M$ in the "middle"
19:     calculate by (5) the polynomials $M_{N-1}$, $M_N$ or $M_{N-1}$, $M_{\omega,N}$ (for c = 1, ω a free parameter)
20:     calculate the N roots $x_i$ of $M_N$ or $M_{\omega,\,N}$ the weights $w_i$ with (6) for each $x_i$
21:     **out:** scale the N 2-tupels ($x_i$, $w_i$) to subinterval [$S_M$]

22: **else**
   comment: 2 subintervals $S_M$ and $S_M + 1$ in the "middle"
23:   calculate by (5) the polynomials $M_{N-1}$, $M_{\omega, N}$ (for c = 0, ω a free parameter) or $M_{N-1}$, $M_N$
24:   calculate the N roots $x_i$ of $M_N$ or $M_{\omega, N}$ ($\alpha_L$, ... , $\alpha_{R, M}$, ..., x) the $\alpha_{R, M}$, ... with (5.1*)
   and the weights $w_i$ with (6) for each $x_i$
25:   **out:** scale the N 2-tupels ($x_i$, $w_i$) to subinterval [$S_M$]
26:   calculate the N roots $x_i$ of $M_N$ or $M_{\omega, N}$ ($\alpha_{L, M}$, …, $\alpha_R$, ... , x) the $\alpha_{L, M}$, … with (5.2*)
   and the weights $w_i$ with (6) for each $x_i$
27:   **out:** scale the N 2-tupels ($x_i$, $w_i$) to subinterval [$S_M + 1$]
28: **end if**
29: **Output:** the set { $x_i$, $w_i$ } $_{i = 1 \ldots NS+1}$ of nodes and weights

_________________________________________________________________

## 9. Some examples for "1/2-rules"

### 9.1. A non uniform, even degree, suboptimal, asymmetric $C^0$ quadrature "1/2-rule"

As in section 5.1 the interval [0, 9] is subdivided in the 6 subintervals [0, 1], [1, 3], [3, 6], [6, 7], [7, 8] and [8, 9], the degree of the rule D = 3

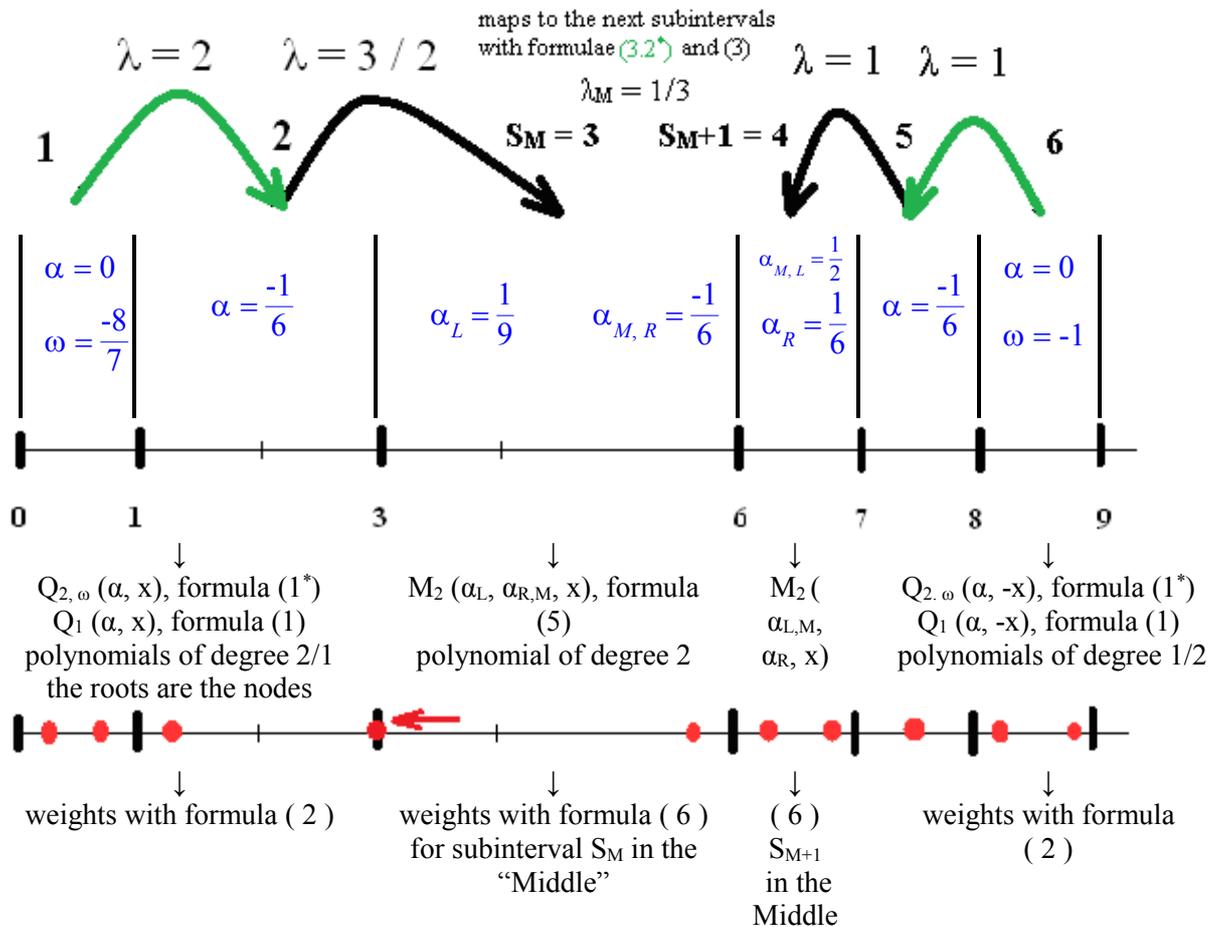

Because the rule is suboptimal, $\alpha_{R,M}$ in subinterval 3 is free and can be chosen so that e.g. a node is at a subinterval boundary e.g. at 3, see red arrow above.

The node/weight list (as rational/algebraic numbers):

(5.1)
```
SQ1 = sqrt (22)
SQ2 = sqrt (113)
SQ3 = sqrt (41)

# subinterval 1
[ 4 / 7 - 1 / 14 * SQ1, - 1 / 44 * SQ1 + 2 / 3 ]
[ 4 / 7 + 1 / 14 * SQ1,   1 / 44 * SQ1 + 2 / 3 ]

# subinterval 2
[ 2, 4 / 3 ]

# subinterval 3
[ 3, 5 / 6 ]   # free parameter chosen, that node at 3.0
[ 9 / 2, 2 ]

# subinterval 4
[181 / 28 - 1 / 28 * SQ2,   3 / 226 * SQ2 + 5 / 6 ]
[181 / 28 + 1 / 28 * SQ2, - 3 / 226 * SQ2 + 5 / 6 ]

# subinterval 5
[ 15 / 2, 2 / 3 ]

# subinterval 6
[ 169 / 20 - 1 / 20 * SQ3,   1 / 164 * SQ3 + 7 / 12 ]
[ 169 / 20 + 1 / 20 * SQ3, - 1 / 164 * SQ3 + 7 / 12 ]
```

## 9.2. A uniform $C^1$ quadrature rule, even degree, $S \to \infty$ i.e. for the real line

The proof is left for the reader again as an exercise.
Hint: There are 2 recursion maps, which map the α, β from subinterval s to the α, β of subinterval s + 2. The interval s can be the subinterval with n or n – 1 nodes.
Taking the subinterval s with n nodes already in the first step when ω with (1.1*) in section 7.1 a square roots is necessary. So the complete recursion map contains a square root, is large and the fixed point is difficult to determine.

So its better to take as first subinterval s the subinterval with n – 1 nodes. To get the desired degree the summand $C_{n-3}$ F (n) / … with lowest index in $Q_{n-1}$ has to be 0, see formula (A.3.1) in appendix A. For F () see formula (1) in section 3. So the α in the first subinterval is a rational function of β. Now apply (3) to get $α_2$, $β_2$ of the second subinterval s + 1. With (3.2*) and an additional independent ω you get $α_3$, $β_3$ of the subinterval s + 2. The recursion map is now rational in $β_1$, ω

Now determine the fixed points of this recursion map, at first determining $β_1$, ω finite, then use formulae (1) and (2). See appendix A.2 with the formula that expresses $Q_n$ (α, β, x) by Gegenbauer polynomials. $Q_n$ (x) has here a factor (1 + x).
There is a second fixed point with ω → ∞, open if this is a quadrature rule, because the limit of one root of $Q_n$ (x) is ∞

The first with finite ω is the optimal quadrature rule give in [14] in section 3.2. formulae (3.2.1) …

## 10. Node distributions in the subintervals

Here I visualize the (combinatorial) node distributions of quadrature rules obtainable with the previous algorithms. The horizontal axis represents the subintervals and the vertical axis the number of nodes in each of the subintervals. It is not claimed that these are all node distributions that occur for $C^0/C^1$ optimal or suboptimal (1-parameter optimal) quadrature rules.

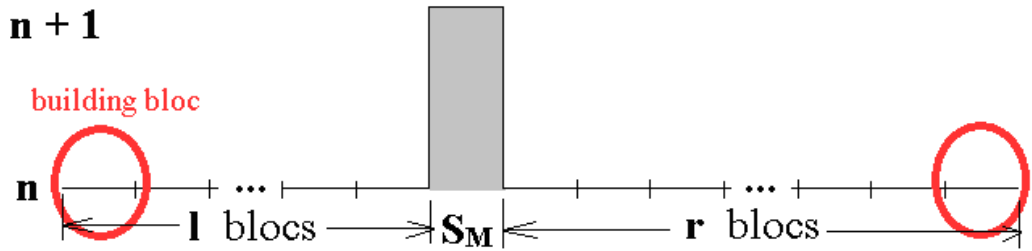

Figure 10.1. $D_1$ (l, r) with S = l + r + 1 subintervals, n or n + 1 nodes per subinterval
a building bloc (in red) consists of 1 subinterval with n nodes

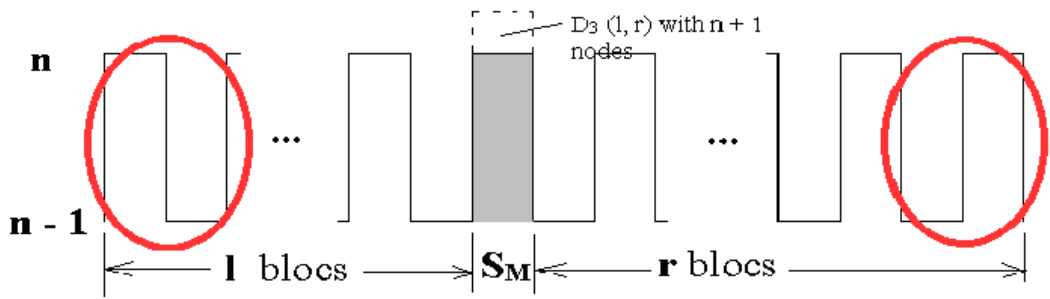

Figure 10.2. $D_2$ (l, r) with S = 2 (l + r) + 1 subintervals, n or n - 1 nodes per subinterval
$D_3$ (l, r) with n + 1 nodes in the subinterval $S_M$
a building bloc (in red) consists of 2 subintervals with n – 1 and n nodes

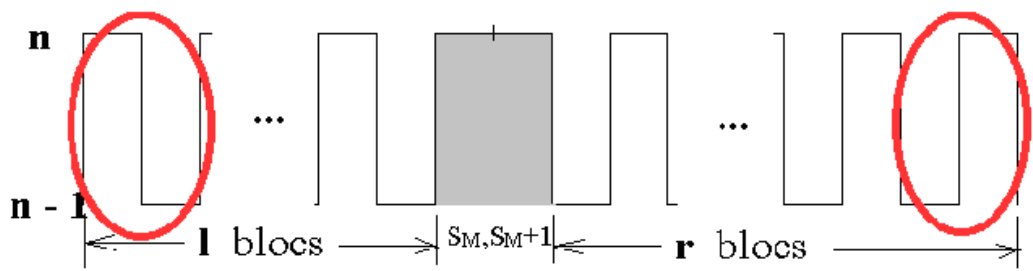

Figure 10.3. $D_4$ (l, r) with S = 2 (l + r) + 2 subintervals, n or n - 1 nodes per subinterval
a building bloc (in red) consists of 2 subintervals with n – 1 and n nodes

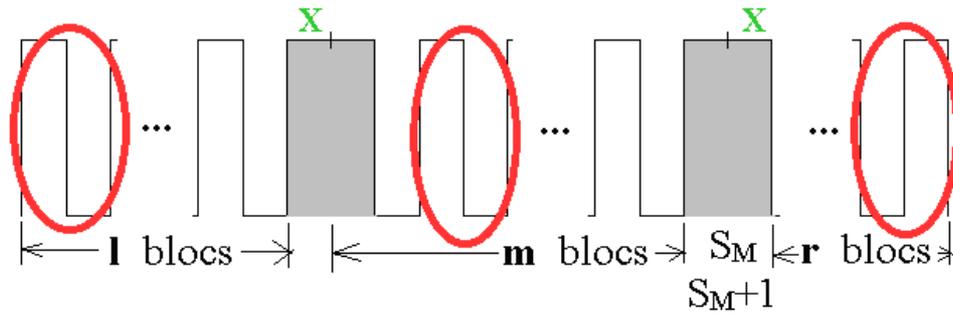

Figure 10.4. $D_5 (l, m, r)$ with $S = 2 (l + m + r) + 3$ subintervals, n or n - 1 nodes per subinterval
a building bloc (in red) consists of 2 subintervals with $n - 1$ and n nodes
The green subintervals X mark positions for a modified "1/2-rules" algorithm
(see 10.2), right of the left X this is a $D_4 (m, r)$ and left of the right X
this is a $D_4 (l, m)$

## 10.1 The 3 node distributions for $C^0$ quadrature rules

| type # | degree D | # of subintervals | optimal no = suboptimal | distribution | remarks |
|---|---|---|---|---|---|
| 1 | even, 2n | even/odd | no | $D_1 (l, r)$ | |
| 2 | | even | no | $D_4 (l, r)$ | |
| 3 | odd, 2n-1 | **odd** | **yes** | **$D_2 (l, r)$** | for non-uniform rules we get different rules for each of the possible $S_M$'s |

Table 10.2.1 The 3 $C^0$ node distributions for optimal and suboptimal rules (optimal in bold)

## 10.2 The 4 node distributions for $C^1$ quadrature rules

| type # | degree D | # of subintervals | optimal no = suboptimal | distribution | remarks |
|---|---|---|---|---|---|
| 1 | odd, 2n+1 | **even/odd** | **yes** | **$D_1 (l, r)$** | |
| 2 | | **even** | **yes** | **$D_4 (l, r)$** | |
| 3 | even, 2n | odd | no | $D_3 (l, r)$ | do not exist for uniform rules, but exist for non-uniform rules under certain conditions |
| 4 | | | no | $D_5 (l, m, r)$ | for $D_5 (0, m, m)$ and $D_5 (0, m, m + 1)$ exist uniform rules for $S \geq 3$ |

Table 10.2.1 The 4 $C^1$ node distributions for optimal and suboptimal rules (optimal in bold)

The algorithm to get quadrature rules of type 4 was not presented until now. The algorithm producing 2 subintervals per step now proceeds form the left or right to the subintervals X,
For this interval we have a free parameter $\omega$ in $Q = Q_n + \omega Q_{n-1}$. $\omega$ is determined e.g. so that Q has a root at -1 or + 1, the $\alpha_N, \beta_N$ for the next segment are given by (3.2*). Because the remaining part is a $D_4 (m, r)$ or $D_4 (l, m)$ the algorithm proceeds as in section 8. Instead of the starting values $\alpha = 0, \beta = 0$ for the algorithm in section 8. we use the $\alpha_N, \beta_N$ from above.

## 11. About the proof, about the verification

The proof uses a CAS just for elementary, symbolic verifications for recursions, integrals of polynomials, neither to solve any system of polynomial or linear equations, nor to do calculations with Groebner bases. All these calculations could be done in finite time by hand (though with a high probability to produce an error). Notice that orthogonal polynomials associated with some semi-classical measures are used. More precisely, they are related with the Beta distribution (the weight function for Jacobi polynomials) as follows:

- H.L. Krall 1940 [7]            with a Dirac δ at only 1 interval end.
- T.H. Koornwinder 1984 [1]      with 2 Dirac δ's at both interval ends.
- J. Arvesú et al 2002 [2]       with additional derivations of Dirac's.

With the definition (2.1) in [1], p. 207, the polynomials in section 2 can be expressed in terms of the Koornwinder polynomials $P_n^{(\alpha, \beta, M, N)}(x)$. Notice that when $\beta = 0$ and $N = 0$, these polynomials were introduced by H. L. Krall.

(10.1)   $M_n(x) = P_n^{(0, 0, \alpha_L, \alpha_R)}(x)$
(10.2)   $Q_n(x) = 1/(n+1)^2 \, P_n^{(1, 0, 2\alpha, 0)}(x)$

According to [2], p. 141, the polynomials in section 3 can be expressed in terms of the orthogonal polynomials $P_n^{(\alpha, \beta, A_1, B_1, A_2, B_2)}(x)$ studied in [2]

(10.3)   $Q_n(x) \sim P_n^{(2, 0, \, 0, \, 8\alpha + 96\beta, \, 0, \, -96\beta)}(x)$
(10.4)   $M_n(x) \sim P_n^{(0, 0, \, 2\alpha, \, 2\alpha, \, 24\beta, \, -24\beta)}(x)$   here given just for the symmetric case

I have verified the results obtained with this formulae numerically with an accuracy of more than 200 digits for several non uniform cases as well as I have checked for all cases the number of subintervals, $S = 1, \ldots, 20$, $N = 1, \ldots, 20$, i.e. degrees up to 40.

I have also compared them with the explicit Gaussian quadrature rules given in [3] for the cubic case $c = 1$, $N = 1$, and in [5] for the quintic case $c = 1$, $N = 2$, as well as with the explicit Gaussian quadrature rules with non uniform knot sequences given in [4] for the cubic case $c = 1$, $N = 1$ (Table 1, geometric stretching, $q = 2$).

*Remark: In [4], Table 1, geometric stretching for N = 6 and 8 seems to be not right, probably due to a copy and paste error, the same first rows as for N = 7 and 9.*

Tests done for the "1/2-rules":

- results of this recipe compared with [13], table 3 p. 18, $N = 3$, subintervals $S = 2, 4, 6, 8, 10$
- results of this recipe compared with [12]
- verified numerically with high accuracy (100 digits) for the uniform and non-uniform cases
  for various N and number of subintervals S

## References

[1]  T.H. Koornwinder. Orthogonal polynomials with weight function
     $(1 - x)^\alpha (1 + x)^\beta + M \delta(x + 1) + N \delta(x - 1)$,
     Canad. Math. Bull. 27 (1984), 205-214
[2]  J. Arvesú, F. Marcellán, R. Álvarez-Nodarse. On a Modification of the
     Jacobi Linear Functional: Asymptotic Properties and Zeros of the
     Corresponding Orthogonal Polynomials
     Acta Applicandae Mathematicae 71: 127-158, 2002

[3]  Geno Nikolov. On certain definite quadrature formulae
     Journal of Computational and Applied Mathematics 75 (1996) 329-343
[4]  Rajid Ait-Haddou, Michael Bartoň, Victor Manuel Calo. Explicit Gaussian quadrature rules for cubic splines with non-uniform knot sequences.
     Preprint, arXiv:1410.7196v1 [math.NA] 27 Oct 2014
[5]  Michael Bartoň, Rajid Ait-Haddou, Victor Manuel Calo. Gaussian quadrature rules for $C^1$ quintic splines.
     Preprint, arXiv:1503.00907v1 [math.NA] 3 Mar 2015
[6]  Michael Bartoň, Victor Manuel Calo. Optimal rules for isogeometric analysis
     Preprint, arXiv:1511.03882v1 [math.NA] 12 Nov 2015
[7]  H.L. Krall. On orthogonal polynomials satisfying a certain fourth order differential equation, The Pennsylvania State College Studies, No. 6, 1940, 1-24
[8]  T.S. Chihara. Orthogonal polynomials and measures with end point masses, Rocky Mountain Journal of Mathematics, Volume 15, Number 3, Summer 1985, 705-719
[9]  K.A. Johannessen. Optimal quadrature for univariate and tensor product splines, Comput. Methods Appl. Mech. Engrg., Volume 316, (2016), 84-99
     http://dx.doi.org/10.1016/j.cma.2016.04.030
[10] K.A. Johannessen. IGA-Quadrature, Open Source, MATLAB sources
     https://github.com/VikingScientist/IGA-quadrature
[11] T. J. R. Hughes, A. Reali, and G. Sangalli. Effcient quadrature for NURBS-based isogeometric analysis. Computer Methods in Applied Mechanics and Engineering, 199 (58): 301-313, 2010.
[12] M. Bartoň and V.M. Calo. Optimal quadrature rules for odd-degree spline spaces and their application to tensor-product-based isogeometric analysis. Computer Methods in Applied Mechanics and Engineering, 305 (2016), 217-240.
[13] Hiemstra, R. R., Calabro, F., Schillinger, D., & Hughes, T. J. (2017). Optimal and reduced quadrature rules for tensor product and hierarchically refined splines in isogeometric analysis. Computer Methods in Applied Mechanics and Engineering, 316, 966-1004.
[14] H. Ruhland. Quadrature rules for $C^0$, $C^1$ splines, the real line and the five (5) families, Preprint, arXiv:1801.03388v4 [math.GM] 10 May 2018

# Appendices

## Appendix A: The $Q_n$ (..., x) and $M_n$ (..., x) expressed by Gegenbauer polynomials

### A.1  c = 0, $Q_n$ (α, x) as sum of 2 Gegenbauer polynomials

With F (n) from (1) section 2.

(A.1.1) $\quad Q_n = \dfrac{C_n\, F(n) + C_{n-1}\, F(n+1)}{n+1}$

### A.2  c = 0, $M_n$ ($α_L$, $α_R$, x) as sum of 3 Gegenbauer polynomials

With H (n) from (5) section 2.

(A.2.1) $\quad M_n = \dfrac{C_n\, H(n) - C_{n-2}\, H(n+1)}{2n+1} + C_{n-1}\, (\alpha_L - \alpha_R)$

### A.3  c = 1, $Q_n$ (α, β, x) as sum of 3 Gegenbauer polynomials

With F (n) from (1) section 3.

$E(n) = 1 + (n+1)(n+2)(\alpha + 3n(n+3)\beta(2 - (n-1)(n+1)(n+2)(n+4)\beta))$

(A.3.1) $\quad Q_n = \dfrac{6\, C_n\, F(n)}{(n+2)(2n+3)} + \dfrac{6\, C_{n-1}\, E(n)}{(n+1)(n+2)} + \dfrac{6\, C_{n-2}\, F(n+1)}{(n+1)(2n+3)}$

### A.4  c = 1, $M_n$ ($α_L$, $β_L$, $α_R$, $β_R$, x) as sum of 5 Gegenbauer polynomials

With H (n) from (5) section 3.

$J_0(n, \alpha, \beta) = 1 + (3 + n + n^2)\alpha + 6(6 + n^2 + 2n^3 + n^4)\beta$
$\qquad - 3(n-3)(n-2)(n-1)n(n+1)(n+2)(n+3)(n+4)\beta^2$

$J_1(n, \alpha, \beta) = 1 + n(n+1)(\alpha + 3(n-1)(n+2)\beta(2 - (n-2)n(n+1)(n+3)\beta))$

$J(n) = \dfrac{1}{2} J_0(n, \alpha_L, \beta_L)\, J_1(n, \alpha_R, \beta_R) + \dfrac{1}{2} J_0(n, \alpha_R, \beta_R)\, J_1(n, \alpha_L, \beta_L)$
$\qquad + 108 (n-1)n(n+1)(n+2)(\beta_L - \beta_R)^2$

$C1 = (\alpha_L - \alpha_R)(3n(\beta_R + \beta_L)(n-1)(n+1)(n+2) - 2)$
$\qquad (3n(\beta_R + \beta_L)(n-2)(n-1)(n+1) - 2)$

$C3 = (\alpha_L - \alpha_R)(3n(\beta_R + \beta_L)(n-1)(n+1)(n+2) - 2)$
$\qquad (3n(\beta_R + \beta_L)(n+1)(n+2)(n+3) - 2)$

$$F13(n) = (\beta_L - \beta_R) n^2 (48 - 144 (n-2)(n+2)(n^2-6)(n-1)^2(n+1)^2 \beta_L \beta_R$$
$$+ 12 (n-1)(n+1)(\alpha_L + \alpha_R) - 48 (n-1)^2(n+1)^2 (\beta_R + \beta_L)$$
$$- 9 (n-2)(n+2)(n-1)^2(n+1)^2 (3 \beta_L \alpha_R + 3 \beta_R \alpha_L + \beta_L \alpha_L + \alpha_R \beta_R))$$

(A.4.1)
$$M_n = \frac{3 C_n H(n)}{(2n+1)(2n+3)} - \frac{6 C_{n-2} J(n)}{(2n-1)(2n+3)} + \frac{3 C_{n-4} H(n+1)}{(2n-1)(2n+1)}$$
$$+ \frac{3}{4} \frac{C_{n-1}(C1 + F13(n)) - C_{n-3}(C3 + F13(n+1))}{2n+1}$$

# Appendix B: The case continuity class c = 0, 1 formulae as Maple code

```
###################### continuity class c = 0 ########################

# formulae (1) : the functions
F (n) = 1 + alpha * n * (n + 1);

Q[n]  = (F (n) + alpha * n) * P[n] + alpha * (1 - x) * D (P[n]);

# formulae (5) : the functions
H (n)    = 1 + n ^ 2 * (alpha[L] + alpha[R] + (n - 1) * (n + 1) * alpha[L] * alpha[R]);
H[1] (n) = alpha[L] + alpha[R] + 2 * n * (n + 1) * alpha[L] * alpha[R];
H[2] (n, alpha) = 1 + alpha * n * (n + 1);

M[n] = (H (n) + n * H[1] (n)) * P[n]
     + (  alpha[L] * H[2] (n, alpha[R]) * (1 - x)
        - alpha[R] * H[2] (n, alpha[L]) * (1 + x)) * D (P[n]);

# formulae (3) :
Gamma = (n + 1) * (1 + n * (n + 2) * alpha);
alpha[new] = (1 + (n + 1) ^ 2 * alpha) / ((n + 1) * Gamma);

#
# additional formulae for the "1/2-rules"
#

# formula (1.1*) :
omega = -  (n * (1 + (n + 1) ^ 2 * alpha) + lambda * (n + 1) *
            (1 + n * (n + 2) * alpha))
          / ((n + 1) * (1 + n ^ 2 * alpha) + lambda * n * (1 + (n - 1) *
            (n + 1) * alpha));

# formulae (3.2*) :
Gamma = (n + 1) * (1 + n * (n + 2) * alpha) + omega * n * (1 + (n - 1) *
        (n + 1) * alpha);
alpha[new] = (n * (1 + (n + 1) ^ 2 * alpha) + omega * (n + 1) * (1 + n ^ 2 * alpha)) /
             (n * (n + 1) * Gamma);

###################### continuity class c = 1 ########################

# formulae (1) : the functions
F (n) = 1 + n * (n + 2) * (alpha + 6 * (n ^ 2 + 2 * n - 1) * beta
                           - 3 * (n - 1) * n * (n + 1) ^ 2 * (n + 2) * (n + 3)
                             * beta ^ 2);
F[1] (n) = 1 * alpha
         + 12 * beta * ((n ^ 2 + 3 * n + 1) - n * (n + 1) ^ 2 * (n + 2) ^ 2
                                              * (n + 3) * beta);
F[2] (n) = beta * (1 - 3 * n * (n + 1) *(n + 2) * (n + 3) * beta);

Q[n] =   (F (n) + n * F[1] (n)) * P[n]  +    F[1] (n) * (1 - x) * D (P[n])
      - 36 * F[2] (n) * D (P[n]) + 12 * F[2] (n) * (1 - x) * D (P[n]) ^ 2;

# formulae (5) : the functions
H[0] (n, alpha, beta) = 1 + n * (n - 1) * (alpha + (n + 1) * (n - 2) * beta
                          * (6 - 3 * beta * (n + 2) * n * (n - 1) * (n - 3)) );
H1L[0] (n) = H[0] (n + 1, alpha[L], beta[L]),
H1R[0] (n) = H[0] (n + 1, alpha[R], beta[R]);
H (n)    = (  H[0] (n, alpha[L], beta[L]) * H1R[0] (n)
            + H[0] (n, alpha[R], beta[R]) * H1L[0] (n) ) / 2
            - 36 * (n - 1) * n ^ 2 * (n + 1) * (beta[L] - beta[R]) ^ 2;

H[1] (n, alpha, beta) = alpha + 12 * n * (n + 1) * beta
                        * (1 - (n - 1) * (n + 2) * (n ^ 2 + n + 3) * beta);
H[2] (n,        beta) = beta * (1 - 3 * (n - 1) * n * (n + 1) * (n + 2) * beta);
H[3] (n, alpha, beta) = H[0] (n + 1, alpha, beta) + 24 * n * (n + 1) * H[2] (n, beta);
H[4] (n, alpha, beta) = 1 + n * (n + 1) * (2 * alpha + 3 * (n - 1) * n * (n + 1)
                          * (n + 2) * (13 * n ^ 2 + 13 * n - 18) * beta ^ 2);
```

```
M[n] = (   H[3] (n, alpha[L], beta[L]) * H1R[0] (n)
         + H[3] (n, alpha[R], beta[R]) * H1L[0] (n)) / 2 * P[n]
         +   (   H[1] (n, alpha[L], beta[L]) * H1R[0] (n) * (1 - x)
               - H[1] (n, alpha[R], beta[R]) * H1L[0] (n) * (1 + x) ) * D (P[n])
       + 12 * (   H[2] (n,            beta[L]) * H1R[0] (n) * (1 - x)
                + H[2] (n,            beta[R]) * H1L[0] (n) * (1 + x) ) * D (P[n]) ^ 2

       # the following part is 0     for beta[L], beta[R] = 0
       #                             or for alpha[L] = alpha[R], beta[L] = beta[R]
       - 36 * (beta[L] - beta[R]) ^ 2 * n * (n + 1)
            * (n * (n + 1) * P[n] - 2 * x * D (P[n]) + 2 * D (P[n]) ^ 2)
       + 12 * (   H[2] (n, beta[R]) * H[4] (n, alpha[L], beta[L])
                - H[2] (n, beta[L]) * H[4] (n, alpha[R], beta[R]) ) * D (P[n])
       + 72 * n * (n + 1) * (beta[L] - beta[R])
            * (beta[L] + beta[R] - 6 * (n - 1) * n * (n + 1) * (n + 2)
                                    * beta[L] * beta[R]) * x * D (P[n]) ^ 2;
```

```
# formulae (3) :
Gamma = (n + 1) * (n + 2) *
           (1 + n * (n + 3) * alpha
             + 6 * n * (n + 3) * (n ^ 2 + 3 * n - 1) * beta
             - 3 * n ^ 2 * (n - 1) * (n + 1) * (n + 2) * (n + 3) ^ 2 * (n + 4) * beta ^ 2
           ) / 2;
E = 1 + (n + 1) * (n + 2) * (  alpha + 3 * n * (n + 3) * beta
                             * (2 - (n - 1) * (n + 1) * (n + 2) * (n + 4) * beta));
G = 1 - 3 * n * (n + 1) * (n + 2) * (n + 3) * beta;

alpha[new] = - alpha + E *
( 4 * (2 * n ^ 2 + 6 * n + 3)
  + n * (n + 3) *
    (
      + (11 * n ^ 2 + 33 * n + 16) * alpha
      + 12 * (4 * n ^ 4 + 24 * n ^ 3 + 34 * n ^ 2 - 6 * n - 8) * beta
      + 3 * n * (n + 1) * (n + 2) * (n + 3) *
        (
          - 4 * (n + 1) * (n + 2) * (2 * n ^ 2 + 6 * n - 5) * beta ^ 2
          - 3 * (n - 1) * n * (n + 1) * (n + 2) * (n + 3) * (n + 4) * alpha * beta ^ 2
          + 2 * (3 * n ^ 2 + 9 * n - 6) * alpha * beta
          + 1 * alpha ^ 2
        )
    )
) / (12 * Gamma ^ 2);

beta[new] = beta + E * G / (6 * (n + 1) * (n + 2) * Gamma);

#
# additional formulae for the "1/2-rules"
#

# formula (1.1*) :
A = alpha^2*(3*lambda^4*n^2*(n-1)^2*(n+2)^2*(n+1)^2+4*lambda^3*n^2*(n-1)*(n+1)
    *(2*n^2+2*n-3)*(n+2)^2+6*lambda^2*n^2*(n-1)*(n+1)*(n^2+2*n-1)*(n+2)^2-n^2
    *(n-1)*(n+3)*(n+2)^2*(n+1)^2)+alpha*(beta^2*(-18*lambda^4*n^3*(n-2)*(n+3)
    *(n-1)^3*(n+2)^3*(n+1)^3-12*lambda^3*n^3*(2*n^2+2*n-9)*(2*n^2+2*n-3)*(n-1)^2
    *(n+1)^2*(n+2)^3-36*lambda^2*n^3*(n^2+n-3)*(n^2+2*n-1)*(n-1)^2*(n+1)^2*(n+2)^3
    +6*n^4*(n+3)*(n-1)^2*(n+2)^3*(n+1)^4)+beta*(36*lambda^4*n^2*(n^2+n-3)*(n-1)^2
    *(n+2)^2*(n+1)^2+24*lambda^3*n^2*(n-1)*(n+1)*(2*n^2+2*n-5)*(2*n^2+2*n-3)*(n+2)^2
    +72*lambda^2*n^2*(n+1)*(n^2+2*n-1)*(n-1)^2*(n+2)^3-12*n^2*(n-1)*(n+3)*(n^2+n-1)
    *(n+2)^2*(n+1)^2)+6*lambda^4*n^2*(n+2)*(n-1)*(n+1)^2+4*lambda^3*n*(n+2)
    *(2*n^2+2*n-3)*(2*n^2+2*n-1)+12*lambda^2*n*(n+2)*(n^2+n-1)*(n^2+2*n-1)-2*n
    *(n-1)*(n+3)*(n+1)*(n+2)^2)+beta^4*(27*lambda^4*n^4*(n-2)^2*(n+3)^2*(n-1)^4*(n+2)^4
    *(n+1)^4+36*lambda^3*n^5*(n-2)*(n+3)*(2*n^2+2*n-9)*(n-1)^3*(n+2)^4*(n+1)^4
    +54*lambda^2*n^5*(n-2)*(n+3)*(n^2+2*n-1)*(n-1)^3*(n+2)^4*(n+1)^4-9*n^6*(n+3)
    *(n-1)^3*(n+2)^4*(n+1)^6)+beta^3*(-108*lambda^4*n^3*(n-2)*(n+3)*(n^2+n-3)
    *(n-1)^3*(n+2)^3*(n+1)^3-144*lambda^3*n^3*(2*n^4+4*n^3-9*n^2-11*n+6)*(n^2+n-3)
    *(n-1)^2*(n+1)^2*(n+2)^3-216*lambda^2*n^3*(n^2+2*n-1)*(n^4+2*n^3-4*n^2-5*n+3)
    *(n-1)^2*(n+1)^2*(n+2)^3+36*n^4*(n+3)*(n^2+2*n-1)*(n-1)^2*(n+2)^3*(n+1)^4)
    +beta^2*(18*lambda^4*n^2*(5*n^4+10*n^3-25*n^2-30*n+54)*(n-1)^2*(n+2)^2*(n+1)^2
    +24*lambda^3*n^2*(n-1)*(n+1)*(10*n^6+30*n^5-35*n^4-120*n^3+67*n^2+132*n-72)*(n+2)^2
    +36*lambda^2*n^2*(n-1)*(n+1)*(n^2+2*n-1)*(5*n^4+10*n^3-15*n^2-20*n+12)*(n+2)^2
    -6*n^2*(n-1)*(n+3)*(5*n^4+10*n^3-5*n^2-10*n+6)*(n+2)^2*(n+1)^2)+beta
    *(36*lambda^4*n^2*(n+2)*(n-1)*(n^2+n-3)*(n+1)^2+48*lambda^3*n*(n-1)
    *(2*n^2+2*n-3)*(n^2+n-1)*(n+2)^2+72*lambda^2*n*(n+2)*(n^2+2*n-1)
    *(n^4+2*n^3-2*n^2-3*n+3)-12*n*(n-1)*(n+3)*(n+1)*(n^2+n-1)*(n+2)^2)+3*lambda^4
    *n^2*(n+1)^2+4*lambda^3*n*(n+2)*(2*n^2+2*n-1)+6*lambda^2*n*(n+2)
    *(n^2+2*n-1)-(n-1)*(n+3)*(n+2)^2;

B = alpha^2*(6*lambda^4*n^2*(n-1)*(n+3)*(n+2)^2*(n+1)^2+8*lambda^3*n^2*(2*n^2+4*n-3)
    *(n+2)^2*(n+1)^2+12*lambda^2*n^2*(n^2+2*n-1)*(n+2)^2*(n+1)^2-2*n^2*(n-1)*(n+3)
    *(n+2)^2*(n+1)^2)+alpha*(beta^2*(-36*lambda^4*n^3*(n-1)^2*(n+3)^2*(n+2)^3*(n+1)^4
    -24*lambda^3*n^3*(n-1)*(n+3)*(4*n^4+16*n^3+20*n^2+8*n-27)*(n+1)^2*(n+2)^3
    -72*lambda^2*n^3*(n-1)*(n+3)*(n^2+2*n+4)*(n^2+2*n-1)*(n+1)^2*(n+2)^3+12*n^3*(n-1)
    *(n+3)*(3+2*n+n^2)*(n+2)^3*(n+1)^4)+beta*(72*lambda^4*n^2*(n-1)*(n+3)*(n^2+2*n-2)
    *(n+2)^2*(n+1)^2+48*lambda^3*n^2*(n+2)^2*(n+1)^2*(2*n^2+4*n-3)^2+144*lambda^2*n^2
    *(n+2)^2*(n+1)^2*(n^2+2*n-1)^2-24*n^3*(n-1)*(n+3)*(n+1)^2*(n+2)^3)
    +12*lambda^4*n*(n+2)*(n^2+2*n-1)*(n+1)^2+8*lambda^3*n*(n+2)*(2*n^2+4*n-1)
    *(2*n^2+4*n+3)+24*lambda^2*n*(n+2)*(n^2+2*n+2)*(n^2+2*n-1)-4*n*(n-1)*(n+3)*(n+2)
    *(n+1)^2)+beta^4*(54*lambda^4*n^4*(n+4)*(n-2)*(n-1)^3*(n+3)^3*(n+2)^4*(n+1)^4
    +72*lambda^3*n^4*(2*n^4+8*n^3-5*n^2-26*n+12)*(n-1)^2*(n+3)^2*(n+2)^4*(n+1)^4
    +108*lambda^2*n^4*(n^2+2*n-2)*(n^2+2*n-1)*(n-1)^2*(n+3)^2*(n+2)^4*(n+1)^4-18*n^5
    *(n-1)^2*(n+3)^2*(n+2)^5*(n+1)^6)+beta^3*(-216*lambda^4*n^3*(n^2+2*n-5)*(n-1)^2
    *(n+3)^2*(n+2)^3*(n+1)^4-288*lambda^3*n^3*(n-1)*(n+3)*(2*n^2+4*n-7)
    *(n^4+4*n^3+4*n^2-3)*(n+1)^2*(n+2)^3-432*lambda^2*n^3*(n-1)*(n+3)*(n^2+2*n+2)
```

```
         *(n^2+2*n-1)*(n^2+2*n-2)*(n+1)^2*(n+2)^3+72*n^3*(n-1)*(n+3)*(n^4+4*n^3+4*n^2-3)
         *(n+2)^3*(n+1)^4)+beta^2*(36*lambda^4*n^2*(n-1)*(n+3)*(5*n^4+20*n^3-13*n^2-66*n+16)
         *(n+2)^2*(n+1)^2+48*lambda^3*n^2*(n-1)*(n+3)*(10*n^4+40*n^3-n^2-82*n+30)*(n+2)^2
         *(n+1)^2+72*lambda^2*n^2*(n^2+2*n-1)*(5*n^4+20*n^3-3*n^2-46*n+36)*(n+2)^2*(n+1)^2
         -12*n^2*(n-1)*(n+3)*(5*n^4+20*n^3+7*n^2-26*n-12)*(n+2)^2*(n+1)^2)
         +beta*(72*lambda^4*n*(n+2)*(n^4+4*n^3+2*n^2-4*n+3)*(n+1)^2+96*lambda^3*n*(n+2)
         *(3+25*n^4-8*n+12*n^5+20*n^3+2*n^6)+144*lambda^2*n*(n+2)*(n^2+2*n-1)
         *(n^4+4*n^3+6*n^2+4*n-2)-24*n*(n-1)*(n+3)*(n+2)*(n+1)^4)+6*lambda^4*n
         *(n+2)*(n+1)^2+8*lambda^3*n*(n+2)*(2*n^2+4*n+3)+12*lambda^2*(n^2+2*n+2)*(n^2+2*n-1)
         -2*n*(n-1)*(n+3)*(n+2);

C  =  alpha^2*(3*lambda^4*n^2*(n+3)^2*(n+2)^2*(n+1)^2+4*lambda^3*n^2*(n+3)*(n+1)
      *(2*n^2+6*n+1)*(n+2)^2+6*lambda^2*n^2*(n+3)*(n+1)*(n^2+2*n-1)*(n+2)^2-n^2*(n-1)
      *(n+3)*(n+2)^2*(n+1)^2)+alpha*(beta^2*(-18*lambda^4*n^3*(n-1)*(n+4)*(n+3)^3
      *(n+2)^3*(n+1)^3-12*lambda^3*n^3*(2*n^2+6*n-5)*(2*n^2+6*n+1)*(n+3)^2*(n+1)^2
      *(n+2)^3-36*lambda^2*n^3*(n^2+2*n-1)*(n^2+3*n-1)*(n+3)^2*(n+1)^2*(n+2)^3+6*n^3
      *(n-1)*(n+3)^2*(n+2)^4*(n+1)^4)+beta*(36*lambda^4*n^2*(n^2+3*n-1)*(n+3)^2*(n+2)^2
      *(n+1)^2+24*lambda^3*n^2*(n+3)*(n+1)*(2*n^2+6*n-1)*(2*n^2+6*n+1)*(n+2)^2
      +72*lambda^2*n^3*(n+1)*(n^2+2*n-1)*(n+3)^2*(n+2)^2-12*n^2*(n-1)*(n+3)*(n^2+3*n+1)
      *(n+2)^2*(n+1)^2)+6*lambda^4*n*(n+3)*(n+2)^2*(n+1)^2+4*lambda^3*n*(n+2)
      *(2*n^2+6*n+1)*(2*n^2+6*n+3)+12*lambda^2*n*(n+2)*(n^2+3*n+1)*(n^2+2*n-1)-2*n^2
      *(n-1)*(n+3)*(n+2)*(n+1))+beta^4*(27*lambda^4*n^4*(n-1)^2*(n+4)^2*(n+3)^4
      *(n+2)^4*(n+1)^4+36*lambda^3*n^4*(n-1)*(n+4)*(2*n^2+6*n-5)*(n+3)^3*(n+1)^4*(n+2)^5
      +54*lambda^2*n^4*(n-1)*(n+4)*(n^2+2*n-1)*(n+3)^3*(n+1)^4*(n+2)^5-9*n^4*(n-1)*(n+3)^3
      *(n+2)^6*(n+1)^6)+beta^3*(-108*lambda^4*n^3*(n-1)*(n+4)*(n^2+3*n-1)*(n+3)^3*(n+2)^3
      *(n+1)^3-144*lambda^3*n^3*(2*n^4+12*n^3+15*n^2-9*n-8)*(n^2+3*n-1)*(n+3)^2*(n+1)^2
      *(n+2)^3-216*lambda^2*n^3*(n^2+2*n-1)*(n^4+6*n^3+8*n^2-3*n-3)*(n+3)^2*(n+1)^2
      *(n+2)^3+36*n^3*(n-1)*(n^2+3*n+1)*(n+3)^2*(n+2)^4*(n+1)^4)+beta^2*(18*lambda^4*n^2
      *(5*n^4+30*n^3+35*n^2-30*n+14)*(n+3)^2*(n+2)^2*(n+1)^2+24*lambda^3*n^2*(n+3)*(n+1)
      *(10*n^6+90*n^5+265*n^4+240*n^3-53*n^2-24*n+12)*(n+2)^2+36*lambda^2*n^2*(n+3)*(n+1)
      *(n^2+2*n-1)*(5*n^4+30*n^3+45*n^2-8)*(n+2)^2-6*n^2*(n-1)*(n+3)
      *(5*n^4+30*n^3+55*n^2+30*n+6)*(n+2)^2*(n+1)^2)+beta*(36*lambda^4*n*(n+3)*(n^2+3*n-1)
      *(n+2)^2*(n+1)^2+48*lambda^3*n^2*(n+3)*(n+2)*(n^2+3*n+1)*(2*n^2+6*n+1)+72*lambda^2*n
      *(n+2)*(n^2+2*n-1)*(n^4+6*n^3+10*n^2+3*n+1)-12*n^2*(n-1)*(n+3)*(n+2)*(n+1)
      *(n^2+3*n+1))+3*lambda^4*(n+2)^2*(n+1)^2+4*lambda^3*n*(n+2)*(2*n^2+6*n+3)
      +6*lambda^2*n*(n+2)*(n^2+2*n-1)-n^2*(n-1)*(n+3);

SQ    = B ^ 2 - 4 * A * C;
omega = (- B + sqrt (SQ)) / 2 / A;
```

```
# formulae (3.2*) :
Gamma = - (n + 2) * (1 + n * (n + 3) * alpha
            + 6 * n * (n + 3) * (n ^ 2 + 3 * n - 1) * beta
            - 3 * n ^ 2 * (n - 1) * (n + 1) * (n + 2) * (n + 3) ^ 2 * (n + 4) * beta ^ 2)
            + omega * (- alpha * n * (n + 2)* (n - 1)
                       + 3 * beta * n * (n + 2) * (n - 1) *
               (beta * n * (n - 2) * (n + 3) * (n + 1) * (n + 2) * (n - 1)
                - 2 * (n ^ 2 + n - 3 )) - n);
A = omega^2*(alpha^2*n*(n-1)*(n+2)*(n+1)*(2*n^2+2*n-3)+alpha*(-3*beta^2*n^2
    *(2*n^2+2*n-9)*(2*n^2+2*n-3)*(n-1)^2*(n+2)^2*(n+1)^2+6*beta*n*(n-2)*(n+1)
    *(2*n^2+2*n-5)*(2*n^2+2*n-3)+(2*n^2+2*n-3)*(2*n^2+2*n-1))
    +9*beta^4*n^4*(n-2)*(n+3)*(2*n^2+2*n-9)*(n+2)^3*(n-1)^3*(n+1)^4
    -36*beta^3*n^2*(2*n^4+4*n^3-9*n^2-11*n+6)*(n^2+n-3)*(n-1)^2*(n+2)^2*(n+1)^2
    +6*beta^2*n*(n-1)*(n+2)*(n+1)*(10*n^6+30*n^5-35*n^4-120*n^3+67*n^2+132*n-72)
    +12*beta*(n+2)*(n-1)*(2*n^2+2*n-3)*(n^2+n-1)+2*n^2+2*n-1)
    +omega*(2*alpha^2*n*(n+2)*(2*n^2+4*n-3)*(n+1)^2
    +alpha*(-6*beta^2*n^2*(n+3)*(n-1)*(4*n^4+16*n^3+20*n^2+8*n-27)*(n+2)^2*(n+1)^2
    +12*beta*n*(n+2)*(n+1)^2*(2*n^2+4*n-3)^2+2*(2*n^2+4*n-1)*(2*n^2+4*n+3))
    +18*beta^4*n^3*(2*n^4+8*n^3-5*n^2-26*n+12)*(n+3)^2*(n-1)^2*(n+2)^3*(n+1)^4
    -72*beta^3*n^2*(n+3)*(n-1)*(2*n^2+4*n-7)*(n^4+4*n^3+4*n^2-3)*(n+2)^2*(n+1)^2
    +12*beta*alpha^2*n*(n-1)*(n+3)*(n+2)*(10*n^4+40*n^3-n^2-82*n+30)*(n+1)^2
    +beta*(72+600*n^4-192*n+288*n^5+480*n^3+48*n^6)+6+8*n+4*n^2)
    +alpha^2*n*(n+3)*(n+2)*(n+1)*(2*n^2+6*n+1)+alpha*(-3*beta^2*n^2*(2*n^2+6*n-5)
    *(2*n^2+6*n+1)*(n+3)^2*(n+2)^2*(n+1)^2+6*beta*n*(n+3)*(n+2)*(n+1)
    *(2*n^2+6*n-1)*(2*n^2+6*n+1)+(2*n^2+6*n+1)*(2*n^2+6*n+3))
    +9*beta^4*n^3*(n-1)*(n+4)*(2*n^2+6*n-5)*(n+3)^3*(n+2)^4*(n+1)^4
    -36*beta^3*n^2*(2*n^4+12*n^3+15*n^2-9*n-8)*(n^2+3*n-1)*(n+3)^2*(n+2)^2*(n+1)^2
    +6*beta^2*n*(n+3)*(n+2)*(n+1)*(10*n^6+90*n^5+265*n^4+240*n^3-53*n^2-24*n+12)
    +12*beta*n*(n+3)*(n^2+3*n+1)*(2*n^2+6*n+1)+2*n^2+6*n+3;

B = omega*(-alpha*n*(n+2)*(n+1)+3*beta^2*n^3*(n-1)*(n+2)^2*(n+1)^3
    -6*beta*n*(n+2)*(n+1)*(n^2+n-1)-2-n)-alpha*n*(n+2)*(n+1)
    +3*beta^2*n^2*(n+3)*(n+2)^3*(n+1)^3-6*beta*n*(n+2)*(n+1)*(n^2+3*n+1)-n;

alpha[new] = 4 / 3 * A / ((n + 1) * Gamma) ^ 2;
beta[new]  = 1 / 3 * B / (Gamma * (n + 1) ^ 2 * (n + 2) * n);

# formulae (5.1*) for 2 subintervals in the middle:

# the 2 defining equations for alpha[R, M] and beta[R, M], see (5.0.1*) .. (5.0.3*)
A (alpha, beta) =
      9 * n ^ 2 * (n + 2) ^ 2 * (n + 1) ^ 4 * (- alpha * (n + 3) * (n - 1)
      + 3 * beta ^ 2 * n * (n + 4) * (n - 2) * (n + 2) * (n + 3) ^ 2 * (n - 1) ^ 2
      - 6 * beta * (n + 3) * (n - 1) * (n ^ 2 + 2 * n - 4) - 1);
B (alpha, beta) =
      18 * n * (n + 2) * (n + 1) ^ 2 * (alpha *(n ^ 2 + 2 * n - 1)
      - 3 * beta ^ 2 * n * (n - 1) * (n + 3) * (n + 2) * (n ^ 2 + 2 * n - 4) * (n + 1) ^ 2
      + 6 * beta * (n ^ 2 + 2 * n - 2) * (n ^ 2 + 2 * n - 1) + 1);
C (alpha, beta) =
      3 * (alpha * (n + 1) ^ 2 - 3 * beta ^ 2 * n ^ 2 * (n + 2) ^ 2 * (n + 1) ^ 4
           + 6 * beta * n * (n + 2) * (n + 1) ^ 2 + 1);
DD (alpha, beta) =
      3 * (n + 1) ^ 2 * (alpha * n * (n + 2)
      - 3 * beta ^ 2 * n ^ 2 * (n + 3) * (n - 1) * (n + 2) ^ 2 * (n + 1) ^ 2
      + 6 * beta * n * (n + 2) * (n ^ 2 + 2 * n - 1) + 1);

A[1]  = A (alpha[L], beta[L]),
B[1]  = B (alpha[L], beta[L]),
C[1]  = C (alpha[L], beta[L]),
DD[1] = DD (alpha[L], beta[L]);

A[2]  = A (alpha[R], beta[R]),
B[2]  = lambda ^ 2 * B (alpha[R], beta[R]),
C[2]  = lambda ^ 4 * C (alpha[R], beta[R]),
DD[2] = - lambda ^ 3 * DD (alpha[R], beta[R]);
```

```
A[3]  = DD[2] * A[1] - DD[1] * A[2],
B[3]  = DD[2] * B[1] - DD[1] * B[2],
C[3]  = DD[2] * C[1] - DD[1] * C[2];

B[4]  = A[2] * B[1]  - A[1] * B[2],
C[4]  = A[2] * C[1]  - A[1] * C[2],
DD[4] = A[2] * DD[1] - A[1] * DD[2];

# formulae (5.1*):
SQ = B[3] ^ 2 - 4 * A[3] * C[3];
beta[R, M]  = (- B[3] - sqrt (SQ)) / 2 / A[3];
alpha[R, M] = - (C[4] + B[4] * beta[R, M]) / DD[4];

alpha[L, M] = - alpha[R, M] / lambda;
beta[L, M]  = + beta[R, M] / lambda ^ 2;
```

# Appendix C: The $Q_n$ (…, x) and $M_n$ (…, x) expressed by Gegenbauer polynomials, the formulae as Maple code

```
# formulae (A.1.1) expressing the Q[n] for c = 0 as sum of 2 Gegenbauer polynomials:

Q[n] = (C[n] * F (n) + C[n - 1] * F (n + 1)) / (n + 1);

# formulae (A.2.1) expressing the M[n] for c = 0 as sum of 3 Gegenbauer polynomials:

M[n] =   (C[n] * H (n) - C[n - 2] * H (n + 1)) / (2 * n + 1)
      + C[n - 1] * (alpha[L] - alpha[R]);

# formulae (A.3.1) expressing the Q[n] for c = 1 as sum of 3 Gegenbauer polynomials:

E (n) = 1 + (n + 1) * (n + 2) * (   alpha + 3 * n * (n + 3) * beta
                                 * (2 - (n - 1) * (n + 1) * (n + 2) * (n + 4) * beta));

Q[n] =   C[n]     * 6 * F (n)     / (n + 2) / (2 * n + 3)
      + C[n - 1] * 6 * E (n)     / (n + 1) / (    n + 2)
      + C[n - 2] * 6 * F (n + 1) / (n + 1) / (2 * n + 3);

# formulae (A.4.1) expressing the M[n] for c = 1 as sum of 5 Gegenbauer polynomials:

J[0] (n, alpha, beta) = 1 + (3 + n + n ^ 2) * alpha + 6 * (6 + n ^ 2 + 2 * n ^ 3
         + n ^ 4) * beta - 3 * (n - 3) * (n - 2) * (n - 1) * n * (n + 1) * (n + 2)
         * (n + 3) * (n + 4) * beta ^ 2;
J[1] (n, alpha, beta) = 1 + n * (n + 1) * (alpha + 3 * (n - 1) * (n + 2) * beta
                       * (2 - (n - 2) * n * (n + 1) * (n + 3) * beta));
J (n) = (   J[0] (n, alpha[L], beta[L]) * J[1] (n, alpha[R], beta[R])
         + J[0] (n, alpha[R], beta[R]) * J[1] (n, alpha[L], beta[L]) ) / 2
      + 108 * (n - 1) * n * (n + 1) * (n + 2) * (beta[L] - beta[R]) ^ 2;

C1 = (alpha[L] - alpha[R]) * (3 * n * (beta[R] + beta[L]) * (n - 1) * (n + 1) * (n + 2) - 2)
                           * (3 * n * (beta[R] + beta[L]) * (n - 2) * (n - 1) * (n + 1) - 2);
C3 = (alpha[L] - alpha[R]) * (3 * n * (beta[R] + beta[L]) * (n - 1) * (n + 1) * (n + 2) - 2)
                           * (3 * n * (beta[R] + beta[L]) * (n + 1) * (n + 2) * (n + 3) - 2);
F13 (n) =   (beta[L] - beta[R]) * n ^ 2 * (48 - 144 * (n - 2) * (n + 2) * (n ^ 2 - 6)
         * (n - 1) ^ 2 * (n + 1) ^ 2 * beta[L] * beta[R] + 12 * (n - 1) * (n + 1)
         * (alpha[L] + alpha[R]) - 48 * (n - 1) ^ 2 * (n + 1) ^ 2 * (beta[R] + beta[L])
         - 9 * (n - 2) * (n + 2) * (n - 1) ^ 2 * (n + 1) ^ 2 * (3 * beta[L] * alpha[R]
         + 3 * alpha[L] * beta[R] + beta[L] * alpha[L] + alpha[R] * beta[R]));

M[n] =   C[n]     * 3 * H (n)     / (2 * n + 1) / (2 * n + 3)
      - C[n - 2] * 6 * J (n)     / (2 * n - 1) / (2 * n + 3)
      + C[n - 4] * 3 * H (n + 1) / (2 * n - 1) / (2 * n + 1)

      + (C[n - 1] * (C1 + F13 (n)) - C[n - 3] * (C3 + F13 (n + 1)))
         * 3 / (2 * n + 1) / 4;
```